\documentclass{scrartcl}
\usepackage{amsmath,amsthm,amscd,amssymb}
\usepackage[mdbch]{mathdesign}
\numberwithin{equation}{section}
\theoremstyle{definition}
\newtheorem{dfn}{Definition}[section]
\newtheorem{example}[dfn]{Example}
\newtheorem{rem}[dfn]{Remark}
\theoremstyle{plain}

\newtheorem{prp}[dfn]{Proposition}
\newtheorem{thm}[dfn]{Theorem}
\newtheorem{cor}[dfn]{Corollary}

\newtheorem{prpdef}[dfn]{Proposition-Definition}

\title{Twisted generalized Weyl algebras, polynomial Cartan matrices and Serre-type relations}
\author{Jonas T. Hartwig}
\date{}

\newcommand\al{\alpha}\newcommand\be{\beta} \newcommand\ga{\gamma}

 \newcommand\la{\lambda}
  \newcommand\si{\sigma}

\newcommand\C{\mathbb{C}}    \newcommand\Z{\mathbb{Z}}   
\newcommand\K{\mathbb{K}}

\DeclareMathOperator{\Aut}{Aut}

\begin{document}
\maketitle
\begin{abstract}
Twisted generalized Weyl algebras (TGWAs) are defined as the quotient of a certain
graded algebra by the maximal graded ideal $I$ with trivial zero
component, analogous to how Kac-Moody algebras can be defined.
In this paper we introduce the class of locally finite TGWAs,
and show that one can associate to such an algebra a
polynomial Cartan matrix (a notion extending the usual
generalized Cartan matrices appearing in Kac-Moody algebra theory)
and that the corresponding generalized 
Serre relations hold in the TGWA.
We also give an explicit construction of a family of locally finite TGWAs
depending on a symmetric generalized Cartan matrix $C$ and some
scalars. The polynomial Cartan matrix of an algebra in this family 
may be regarded as a deformation of the original matrix $C$ and
gives rise to quantum Serre relations in the TGWA.
We conjecture that these relations generate the graded ideal $I$ for these
algebras, and prove it in type $A_2$.
\end{abstract}

{\bf Mathematics Subject Classification (MSC2010):} 16S32 (17B37) \\ 


\section{Introduction}
Generalized Weyl algebras (GWAs) were introduced by Bavula \cite{B},
and by Rosenberg \cite{R} under the name hyperbolic rings.
Their structure and representation theory have been extensively
studied in varying degrees of generality,
see \cite{BJ},\cite{BvO},\cite{BB},\cite{BBF},\cite{DGO},\cite{H08}
and references therein.
Examples of generalized Weyl algebras include
the $n$:th Weyl algebra $A_n$, the enveloping algebra $U(\mathfrak{sl}_2)$
and the quantum group $U_q(\mathfrak{sl}_2)$ as well as many other
interesting families of algebras
(see for example \cite{CL} and references therein).

In an attempt to enlarge
the class of GWAs so as to include also the enveloping algebras of semisimple
Lie algebras of higher rank,
Mazorchuk and Turowska \cite{MT} defined the notion of a
twisted generalized Weyl algebra (TGWA).
They are natural generalizations
of the GWAs but their structure is more complicated.
Representations of TGWAs were investigated in
 \cite{MT},\cite{MPT},\cite{MT02} and \cite{H}.

In \cite{MT} a first indication of a relation to higher rank Lie algebras
was found (in the support of weight modules).
Later, in \cite{MPT}, it was shown that the
Mickelsson step algebras $Z(\mathfrak{gl}_{n+1},\mathfrak{gl}_n\oplus \mathfrak{gl}_1)$
as well as some extended orthogonal Gelfand-Zetlin (OGZ) algebras 
associated to $\mathfrak{gl}_n$ are examples of TGWAs.
It is still unknown whether enveloping algebras of higher rank Lie algebras
are examples of TGWAs, but it is known that a certain
 localization of $U(\mathfrak{gl}_n)$ is isomorphic to an extended OGZ algebra,
 and hence is a TGWA.

In this paper we further strengthen the link between TGWAs and
enveloping algebras of semisimple Lie algebras.
Namely, we define a natural subclass of TGWAs which we call 
\emph{locally finite} TGWAs (Definition \ref{dfn:locfin}). We show
(Theorem \ref{thm:lfthm}) that in
those algebras, relations of the following type hold:
\[
 X_i^{m_{ij}}X_j+\eta_{ij}^{(1)} X_i^{m_{ij}-1}X_jX_i+
 \cdots + \eta_{ij}^{(m_{ij})} X_jX_i^{m_{ij}}=0,
\]
where $\eta_{ij}^{(k)}$ 
are scalars. Such relations may be regarded as generalized Serre relations.
Furthermore, we construct examples
 of locally finite TGWAs, denoted $\mathcal{T}_{q,\mu}(C)$,
were these relations are the quantum Serre relations associated
to a symmetric generalized Cartan matrix $C$.
More precisely, Theorem \ref{thm:TqmuC}c) says that
there exist nonzero algebra homomorphisms
\[\varphi_\pm : U_q(\mathfrak{n}_\pm)\to \mathcal{T}_{q,\mu}(C)\]
where $\mathfrak{n}_\pm$ are the positive and negative part of
the Kac-Moody algebra $\mathfrak{g}(C)$ associated to $C$.
We conjecture that the maps $\varphi_\pm$ are injective and prove it
in the case when $C$ is of type $A_2$ (Corollary \ref{cor:A2}).
We believe that the algebras in this family merit further study
and we plan to investigate their structure and representations
in forthcoming papers.

The plan of this paper is as follows.
In Section \ref{sec:notation} we recall the definition of TGWAs and
give some examples.
After establishing some preliminary results in Section \ref{sec:preliminary}
we define the class of locally finite TGWAs in Section \ref{sec:locallyfinite} and
show how to associate polynomial Cartan matrices to them such that
the corresponding Serre relations hold.
Then, in Section \ref{sec:cartan}, we give the construction
of the locally finite TGWAs, $\mathcal{T}_{q,\mu}(C)$,
associated to a symmetric generalized Cartan matrix and some parameters
 $q$, $\mu=(\mu_{ij})$.
Finally, in Section \ref{sec:caseA2} we give sufficient conditions for
the Serre-type relations to generate the maximal graded ideal appearing in
the definition of TGWAs, which in particular can be used to obtain
an explicit presentation of the
algebra $\mathcal{T}_{q,\mu}(C)$ when $C$ is of type $A_2$.

\section*{Acknowledgements}
The author was supported by the
Netherlands Organization for Scientific Research (NWO) in the
VIDI-project ``Symmetry and modularity in exactly solvable models''.
The author would like to thank L. Turowska for commenting on an early version
of this paper, and J. Palmkvist and J. \"Oinert for interesting discussions.

\section{Notation and definitions} \label{sec:notation}

We fix an arbitrary ground field $\K$. All algebras are associative
unital $\K$-algebras.
Let us recall the definition of a twisted generalized Weyl algebra \cite{MT,MPT}.
The input for this construction is
a positive integer $n$, a commutative $\K$-algebra $R$,
an $n$-tuple $\si=(\si_1,\ldots,\si_n)$
of commuting $\K$-algebra automorphisms of $R$,
an $n$-tuple $t=(t_1,\ldots,t_n)\in (R\backslash\{0\})^n$,
and a symmetric matrix $\mu=(\mu_{ij})_{i,j=1}^n$ with entries from
$\K\backslash\{0\}$ but diagonal elements $\mu_{ii}$ left undefined (they are irrelevant).
The following consistency relation is required:
\begin{equation}\label{eq:consistency}
t_it_j=\mu_{ij}\mu_{ji}\si_j^{-1}(t_i)\si_i^{-1}(t_j)\qquad\forall i\neq j.
\end{equation}
The associated \emph{twisted generalized Weyl construction (TGWC)},
$A'=A'(R,\si,t,\mu)$,
is the algebra obtained from $R$ by adjoining new
noncommuting generators $X_1,\ldots,X_n$, $Y_1,\ldots,Y_n$
with the following defining relations
for $i,j=1,\ldots,n$, $i\neq j$:
\begin{subequations}\label{eq:tgwcrels}
\begin{align}
\label{eq:tgwarels1}
X_ir  &=\si_i(r)X_i,  & Y_ir&=\si_i^{-1}(r)Y_i, 
 && \text{ $\forall r\in R$,} \\
\label{eq:tgwarels2}
Y_iX_i&=t_i, & X_iY_i&=\si_i(t_i),
 && \\
\label{eq:tgwarels3}
X_iY_j&=\mu_{ij}Y_jX_i. &&
 &&
\end{align}
\end{subequations}
One can show that relation \eqref{eq:consistency} is sufficient
for the algebra $A'$ to be nontrivial, and also necessary if the $t_i$
are not zero-divisors in $R$.
$A'$ has a $\Z^n$-gradation $\{A'_g\}_{g\in\Z^n}$
given by $\deg X_i=e_i, \deg Y_i=-e_i, \deg r=0 \forall r\in R$,
where $\{e_i\}_{i=1}^n$ is the standard $\Z$-basis in $\Z^n$.
One can check that $A'_0=R$.
Let $I\subseteq A'$ be the unique maximal graded ideal intersecting $R$ trivially.
Equivalently $I$ is the sum of all graded ideals intersecting $R$ trivially.
The \emph{twisted generalized Weyl algebra (TGWA)},
 $A=A(R,\si,t,\mu)$, is defined as the
quotient $A=A'/I$. Since $I$ is graded,
$A$ inherits a $\Z^n$-gradation from $A'$.
The images of the elements $X_i,Y_i$ in $A$ will be denoted by the same letters.

In this paper, the assumption that $\mu$ is symmetric is a matter of
convenience, because then relations \eqref{eq:tgwcrels} imply that $A'$
carries an anti-involution $\ast$, (that is, a $\K$-linear map $A'\to A'$
 with $a^{\ast\ast}=a$, $(ab)^\ast=b^\ast a^\ast$)
determined by
\begin{equation}\label{eq:involution}
X_i^\ast=Y_i,\quad Y_i^\ast=X_i,\quad i=1,\ldots,n,\qquad r^\ast=r\;\; \forall r\in R.
\end{equation}
Clearly $I^\ast=I$ so $\ast$ descends to an anti-involution on $A$.

By a \emph{homogenous} element in $A'$ (resp. $A$) we mean
an element of $\cup_{g\in\Z^n}A'_g$ (resp. $\cup_{g\in\Z^n}A_g$).
By a \emph{monic monomial} in $A'$ (or $A$) we mean a product
$Z_1\cdots Z_k$ where $Z_j\in \{X_i\}_{i=1}^n\cup\{Y_i\}_{i=1}^n$
for each $j$.
The group $\Z^n$ acts on $R$ via the automorphisms $\si_i$:
$g(r)=(\si_1^{g_1}\cdots\si_n^{g_n})(r)$ for $g=(g_1,\ldots,g_n)
\in\Z^n$ and $r\in R$. Using this action and \eqref{eq:tgwarels1}
we have
\begin{equation}
a\cdot r=(\deg a)(r)\cdot a
\quad \text{for any homogenous $a\in A$ and any $r\in R$}.
\end{equation}

\begin{example}[The TGWA of ``type $A_2$'' from Example 3 in \cite{MT}]
\label{ex:typeA2}
Let $n=2$, $R=\K[H]$, $\si_1(H)=H+1$, $\si_2(H)=H-1$,
and $t_1=H$, $t_2=H+1$, and $\mu_{12}=\mu_{21}=1$.
Then the consistency relation $t_it_j=\si_j^{-1}(t_i)\si_i^{-1}(t_j)$, $i\neq j$ hold.
Let $A=A(R,\si,t,\mu)$ be the corresponding TGWA. 
In \cite{MT} this example was associated to the coxeter graph of type $A_2$,
but the ideal $I$ was not described.
We will come back to this example and eventually
exhibit a set of generators
for the ideal $I$  (Example \ref{ex:typeA2serre}
and Example \ref{ex:typeA2generators}).
\end{example}

\begin{example}[Quantized Weyl algebras]
Let $\bar q=(q_1,\ldots,q_n)$ be an $n$-tuple of elements
of $\K\backslash\{0,1\}$.
Let $\Lambda=(\la_{ij})_{i,j=1}^n$ be an $n\times n$ matrix with
nonzero entries from $\K$ such that $\la_{ji}=\la_{ij}^{-1}$ for $i<j$.
The $n$:th quantized Weyl algebra
$A_n^{\bar q,\Lambda}$ is the $\K$-algebra with generators
$x_i,y_i$, $1\le i\le n$, and relations
\begin{align}
x_ix_j&=q_i\la_{ij}x_jx_i, & y_iy_j&=\la_{ij}y_jy_i,
\label{eq:qweylrel1}\\
x_iy_j&=\la_{ji}y_jx_i, & x_jy_i&=q_i\la_{ij}y_ix_j,
\label{eq:qweylrel2}
\end{align}
for $1\le i<j\le n$, and
\begin{equation}
\label{eq:qweylrel3}
x_iy_i-q_iy_ix_i=1+\sum_{k=1}^{i-1}(q_k-1)y_ix_i,
\end{equation}
for $i=1,\ldots,n$.
This algebra was introduced in \cite{PW} (for special parameters)
and investigated since then by many authors.
Its representation theory has been studied from the point of view of TGWAs in
\cite{MT02} and \cite{H}.
To realize it as a TGWA,
let $R=\K[t_1,\ldots, t_n]$ be the polynomial algebra in
$n$ variables and $\si_i$ the $\K$-algebra
automorphisms of $R$ defined by
\begin{equation}\label{eq:qweylsigmadef}
\si_i(t_j)=
\begin{cases}
t_j, & j<i, \\
1+q_it_i+\sum_{k=1}^{i-1}(q_k-1)t_k, & j=i, \\
q_it_j, & j>i.
\end{cases}
\end{equation}
One can check that the $\si_i$ commute.
Let $\mu=(\mu_{ij})_{i,j=1}^n$
be defined by $\mu_{ij}=\la_{ji}$ and $\mu_{ji}=q_i\la_{ij}$
for $i<j$ (so if the involution \eqref{eq:involution} is required,
one needs to impose also that $\la_{ji}=q_i\la_{ij}$, $i<j$). Put
$\si=(\si_1,\ldots,\si_n)$ and
$t=(t_1,\ldots,t_n)$.
One can show that the maximal graded ideal of
the TGWC
$A'(R,\si,t,\mu)$
is generated by the elements
\[X_iX_j-q_i\la_{ij}X_jX_i,\quad Y_iY_j-\la_{ij}Y_jY_i,\quad
1\le i<j\le n\]
so that $A_n^{\bar q,\Lambda}$ is isomorphic
to the TGWA $A(R,\si,t,\mu)$
via $x_i\mapsto X_i$, $y_i\mapsto Y_i$.
\end{example}

\section{Preliminary results} \label{sec:preliminary}
The following observation about the ideal $I$ in a TGWC is very useful.
Recall that an ideal $J$ in a graded algebra $B$ is called
\emph{completely gr-prime} if $ab\in J$ implies $a\in J$
or $b\in J$ for any homogenous elements $a,b\in B$.
\begin{prp} \label{prp:completely_gr-prime}
Let $A'=A'(R,\si,t,\mu)$ be a TGWC,
where $R$ is a domain. Then $I$ is completely gr-prime.
\end{prp}
\begin{proof} 
Let $a,b\in A'$ be homogenous and suppose $ab\in I$.
First we assume that $a=b^\ast$.
Then, since $\deg(b^\ast)=-\deg(b)$, we have $b^\ast b\in I\cap R=0$.
Let $c_1,c_2\in A'$ be any homogenous elements such
that $\deg(c_1 b c_2)=0$. Then
$(c_1bc_2)^\ast\cdot (c_1bc_2) = (\deg c_2^\ast b^\ast) (c_1^\ast c_1)
\cdot (\deg c_2^\ast)(b^\ast b) \cdot c_2^\ast c_2 = 0$.
Since $R$ has no zero-divisors we conclude that
$c_1bc_2=0$. This proves that $b\in I$.

Now let $a,b$ be any homogenous element in $A'$ with $ab\in I$.
Then $0=(ab)^\ast\cdot ab=
b^\ast (a^\ast a) b = (\deg b^\ast)(a^\ast a) b^\ast b$.
So since $R$ has no zero-divisors either $a^\ast a=0$
or $b^\ast b=0$. Thus, by the first part, either $a\in I$ or $b\in I$.
\end{proof}

\begin{rem}
Proposition \ref{prp:completely_gr-prime}
can be viewed as a refinement of the statement that
the \emph{Shapovalov form}
$F:A\times A\to R$, given by $F(a,b)=\mathfrak{p}(a^\ast b)$,
where $\mathfrak{p}: A\to A_0=R$ is the graded projection,
is non-degenerate when $A=A(R,\si,t,\mu)$ is a TGWA
in which $R$ is a domain,
see \cite[Corollary 5.5]{MPT}.
\end{rem}

Since the $t_i$ are nonzero, we get the following corollary.
\begin{cor}\label{cor:nomonicmonomials}
Let $A'=A'(R,\si,t,\mu)$ be a TGWC,
where $R$ is a domain. Then no monic monomial belongs to $I$.
\end{cor}

Using Proposition \ref{prp:completely_gr-prime}
we can prove the following explicit correspondence
between certain relations in $R$ and certain elements in the ideal $I$ of a TGWC.
\begin{prp}\label{prp:relequiv}
Let $A'=A'(R,\si,t,\mu)$ be a TGWC and
suppose $R$ is a domain. Fix $i,j\in\{1,\ldots,n\}$, $i\neq j$.
Let $m\in\Z_{> 0}$ and $r_0,\ldots,r_m\in R$. Then the following
are equivalent:
\begin{enumerate}
\item[(i)] $r_0 X_i^m X_j + r_1 X_i^{m-1}X_jX_i + \cdots + r_m X_iX_j^m \in I$,
\item[(ii)] $Y_jY_i^m r_0 + Y_iY_jY_i^{m-1} r_1 + \cdots Y_i^m Y_j r_m \in I$,
\item[(iii)] $s_0\si_i^m(t_j) + s_i \si_i^{m-1}(t_j)+\cdots + s_m t_j=0$,
\qquad where $s_k=\mu_{ij}^k\si_j^{-1}(r_k)$ for $k=0,\ldots,m$.
\end{enumerate}
\end{prp}
\begin{proof}
By applying the involution $\ast$, (i) is clearly equivalent to (ii).
Put $a=\sum_{k=0}^m r_k X_i^{m-k}X_jX_i^k$. Then $a$ is homogenous
of degree $me_i+e_j$, where $e_i=(0,\ldots,1,\ldots,0)$ ($1$ on $i$:th position).
So if $a\in I$ then $aY_jY_i^m=0$ by definition of $I$.
Conversely, if $aY_jY_i^m=0$, then $a\in I$
since $I$ is completely gr-prime and does not contain the monic monomial $Y_jY_i^k$.
On the other hand,
\begin{align*}
aY_jY_i^m&=\Big(\sum_{k=0}^m r_k X_i^{m-k}X_jX_i^k\Big)\cdot Y_jY_i^m
= \Big(\sum_{k=0}^m \mu_{ij}^k r_k \si_i^{m-k}\big(\si_j(t_j)\big)\Big)
 \cdot X_i^mY_i^m=\\
&= \si_j\Big(\sum_{k=0}^m \mu_{ij}^k\si_j^{-1}(r_k)\cdot \si_i^{m-k}(t_j)
 \Big)\cdot X_i^mY_i^m,
\end{align*}
which is zero iff (iii) holds, since $X_i^mY_i^m\in R$ and is nonzero,
and $R$ has no zero-divisors.
\end{proof}

\subsection{A note on the case of Noetherian base ring $R$}
This subsection is independent of the rest of the paper.
We include it to show that one can say something about the ideal
$I$ even under very mild conditions on the data $R,\si,t,\mu$.

\begin{prp}
Let $A'=A'(R,\si,t,\mu)$ be a TGWC, and
assume that $R$ is a Noetherian domain. Then
for any $i,j\in\{1,\ldots,n\}$, $i\neq j$, there exist
 $m_{ij},n_{ij}\in\Z_{>0}$
and $r_{ij}^{(1)},\ldots,r_{ij}^{(m_{ij})}\in R$ and $s_{ij}^{(1)},\ldots 
s_{ij}^{(n_{ij})}\in R$, where $r_{ij}^{(m_{ij})},s_{ij}^{(n_{ij})}\neq 0$,
such that
\begin{equation}\label{eq:genserre1}
X_i^{m_{ij}}X_j + r_{ij}^{(1)} X_i^{m_{ij}-1} X_j X_i + \cdots +
 r_{ij}^{(m_{ij})} X_jX_i^{m_{ij}}\in I
\end{equation}
and
\begin{equation}\label{eq:genserre2}
X_jX_i^{n_{ij}} + s_{ij}^{(1)} X_i X_j X_i^{n_{ij}-1} + \cdots +
 s_{ij}^{(n_{ij})} X_i^{n_{ij}}X_j\in I.
\end{equation}
If $m_{ij}$ and $n_{ij}$ are assumed to be minimal, and if all
 $r_{ij}^{(k)}$ ($k=1,\ldots,m_{ij}$) are invertible
  or if all $s_{ij}^{(k)}$ ($k=1,\ldots, n_{ij}$) are invertible,
then the elements in \eqref{eq:genserre1} and \eqref{eq:genserre2}
are multiples of each other by invertible elements from $R$, and
the coefficients are uniquely determined.
By applying the anti-involution $\ast$ we have analogous identities
for the $Y_i$'s.
\end{prp}
\begin{rem}
Note that the leftmost coefficients in \eqref{eq:genserre1},\eqref{eq:genserre2}
are $1$. Otherwise the
statement would be trivial due to the identity
$\si_j(t_j)X_iX_j=\mu_{ij}^{-1} X_jX_i t_i$ which follows from
\eqref{eq:tgwcrels}.
\end{rem}
\begin{proof}
Since $R$ is Noetherian, the ascending chain of ideals
\[(t_j)\subseteq \big(t_j,\si_i(t_j)\big)\subseteq 
\big(t_j,\si_i(t_j),\si_i^2(t_j)\big)\subseteq \cdots \]
in $R$ must stabilize. Thus there is a positive integer $k$
and $r_1',\ldots,r_k'\in R$ such that
\begin{equation}\label{eq:serreproof1}
\si_i^k(t_j)= r_1'\si_i^{k-1}(t_j) +r_2'\si_i^{k-2}(t_j)+\cdots
 +r_k't_j.
\end{equation}
If $k$ was chosen minimal, then $r_k'\neq 0$. Indeed, otherwise, apply
$\si_i^{-1}$ to \eqref{eq:serreproof1} to contradict minimality.
By Proposition \ref{prp:relequiv}, we get a relation of the form
\eqref{eq:genserre1} with minimal $m_{ij}$.
For \eqref{eq:genserre1}, look instead on the chain
$(t_j)\subseteq (t_j,\si_i^{-1}(t_j))\subseteq \cdots $
and apply a positive power of $\si$.

Assume $m_{ij}$ and $n_{ij}$ are minimal and, say, all $r_{ij}^{(k)}$
are invertible. Suppose we have another identity of the form
\eqref{eq:genserre1} with the same $m_{ij}$. Then we can eliminate the
last term and get that $X_i$ multiplied by a certain element $b$
of degree $(m_{ij}-1)e_i+e_j$ lies in $I$,
which by Proposition \ref{prp:completely_gr-prime} gives that $b\in I$,
contradicting minimality of $m_{ij}$.
Then by uniqueness \eqref{eq:genserre1}
and \eqref{eq:genserre2} are just multiples of eachother.
\end{proof}

\section{Locally finite TGWAs and their polynomial Cartan matrices}
\label{sec:locallyfinite}

In view of Proposition \ref{prp:relequiv},
the following class of TGWAs seems natural to consider.
\begin{dfn}\label{dfn:locfin}
\begin{enumerate}
\item[a)] Let $A=A(R,\si,t,\mu)$ be a TGWA.
Define the following $\K$-linear subspaces of $R$ for $i,j=1,\ldots,n$:
\begin{equation}\label{eq:Vijdef}
V_{ij} = \mathrm{Span}_{\K}
 \{\si_i^k(t_j)\mid  k\in\Z\}.
\end{equation}
We say that $A$ is
\emph{locally finite} if $\dim_\K V_{ij}<\infty$ for all $i,j$.
\item[b)]
 To a locally finite TGWA, $A$, we associate the
\emph{matrix of minimal polynomials}
 $P_A=(p_{ij})_{i,j=1}^n$ by
letting $p_{ij}\in\K[x]$ be the minimal polynomial for $\si_i$ acting
on the finite-dimensional space $V_{ij}$. 
Equivalently, $p_{ij}$ is the monic polynomial of minimal degree such that
$\big(p_{ij}(\si_i)\big)(t_j)=0$.
\end{enumerate}
\end{dfn}
\begin{rem}
If the algebra $R$ is generated, as a $\K$-algebra, by the 
set $\{\si_i^k(t_j)\mid i,j=1,\ldots,n,\; k\in\Z\}$, then locally finiteness
of the TGWA is equivalent to that the automorphisms $\si_i$
act locally finitely on $R$, hence the choice of terminology.
\end{rem}

Recall that a generalized Cartan matrix is a matrix $C=(c_{ij})_{i,j=1}^n$
with $c_{ij}\in \Z\,\forall i,j$ which satisfies
i) $c_{ii}=2\,\forall i$, ii) $c_{ij}\le 0\,\forall i\neq j$,
 iii) $c_{ij}=0 \Leftrightarrow c_{ji}=0\,\forall i\neq j$.
We make the following definition.
\begin{dfn} A \emph{polynomial Cartan matrix},
 $P=(p_{ij})_{i,j=1}^n$,
 is a matrix of polynomials $p_{ij}\in\K[x]$ such that
 its shifted degree matrix $(1-\deg p_{ij})_{i,j=1}^n$ coincides
 with a generalized Cartan matrix away from the diagonal. Equivalently, $P$
 is a polynomial Cartan matrix iff for any $i\neq j$ we have
 $\deg p_{ij}\ge 1$ and $\deg p_{ij}=1$ iff $\deg p_{ji}=1$.
 We will denote the Cartan matrix obtained from $P$ by $C(P)$.
\end{dfn}

The point in making these definitions is the following theorem.

\begin{thm}
\label{thm:lfthm}
Let $A=A(R,\si,t,\mu)$ be a locally finite TGWA, where $R$ is a
domain. Then
\begin{enumerate}
\item[a)] the matrix of minimal polynomials $P_A=(p_{ij})_{i,j=1}^n$
 of $A$ is a polynomial Cartan matrix,
\item[b)] writing
\[p_{ij}(x)=x^{m_{ij}}+\la_{ij}^{(1)}x^{m_{ij}-1}+\cdots+\la_{ij}^{(m_{ij})},\]
where all $\la_{ij}^{(k)}\in\K$, the following identities hold in $A$, for any $i\neq j$:
\begin{equation}\label{eq:genserre3}
 X_i^{m_{ij}}X_j+\la_{ij}^{(1)}\mu_{ij}^{-1} X_i^{m_{ij}-1}X_jX_i+
 \cdots + \la_{ij}^{(m_{ij})} \mu_{ij}^{-m_{ij}}X_jX_i^{m_{ij}}=0
\end{equation}
and
\begin{equation}\label{eq:genserre3Y}
 Y_jY_i^{m_{ij}}+\la_{ij}^{(1)}\mu_{ij}^{-1} Y_iY_jY_i^{m_{ij}-1}+
 \cdots + \la_{ij}^{(m_{ij})} \mu_{ij}^{-m_{ij}}Y_i^{m_{ij}}Y_j=0.
\end{equation}
Moreover, for any $i\neq j$ and $m<m_{ij}$, the sets
$\{X_i^{m-k}X_jX_i^k\}_{k=0}^m$ and
$\{Y_i^{m-k}Y_jY_i^k\}_{k=0}^m$ are linearly independent over $\K$.
\end{enumerate}
\end{thm}
\begin{proof}
a) Since $t_i\neq 0$ for each $i$ it is clear that $\deg p_{ij}\ge 1$
for all $i,j$.  And if $p_{ij}$ has degree
one for some $i\neq j$, then $\si_i(t_j)=\la t_j$ for some $\la\in\K\backslash\{0\}$.
Then the consistency relation \eqref{eq:consistency} and that
$t_j$ is not a zero-divisor imply that $\si_j(t_i)=\la't_i$
for some $\la'\in\K\backslash\{0\}$. Thus $\deg(p_{ji})=1$.

b) This is immediate from Proposition \ref{prp:relequiv}.
\end{proof}

Note that the constant terms $\la_{ij}^{(m_{ij})}$ are all nonzero, since
the $p_{ij}$ are minimal polynomials for invertible linear transformations.

\begin{example}[The TGWA of ``type $A_2$'', contd.]
\label{ex:typeA2serre}
Recall that $n=2$, $R=\K[H]$, $\si_1(H)=H+1$, $\si_2(H)=H-1$,
and $t_1=H$, $t_2=H+1$, and $\mu_{12}=\mu_{21}=1$ and
$A=A(R,\si,t,\mu)$ is the corresponding TGWA.
Clearly $A$ is locally finite with $V_{ij}=\C H\oplus \C 1$ for $i,j=1,2$.
Observing that $\si_2(t_1)$ and $t_1$ are linearly independent and that
\[\si_2^2(t_1)-2\si_2(t_1)+t_1=H-2 - 2(H-1) + H = 0\]
we see that the minimal polynomial $p_{21}$ for $\si_2$ acting on $V_{21}$
is given by $p_{21}(x)=x^2-2x+1=(x-1)^2$. Similarly one checks that in fact
$p_{ij}(x)=(x-1)^2$ for all $i,j=1,2$. Thus, by Theorem \ref{thm:lfthm}b),
we have the relations
\[X_1^2X_2-2X_1X_2X_2+X_2X_1^2=0,\quad Y_1^2Y_2-2Y_1Y_2Y_2+Y_2Y_1^2=0,\]
in $A$, which are precisely the Serre relations in the enveloping algebra
of $\mathfrak{sl}_3(\K)$. Also,
 $X_1X_2$, $X_2X_1$, and $Y_1Y_2$, $Y_2Y_1$ are linearly independent.
Note that $C(P_A)=\begin{bmatrix}2&-1\\-1&2\end{bmatrix}$ is the
Cartan matrix of type $A_2$, thus really motivating us to call this algebra
a TGWA of type $A_2$.
\end{example}
\begin{example}
From \eqref{eq:qweylsigmadef} it is easy to see that
the quantized Weyl algebra
 $A_n^{\bar q,\Lambda}$ is also locally finite, and that
the Cartan matrix $C(P)$ associated to its matrix of minimal polynomials $P$
is of the type $(A_1)^n=A_1\times\cdots\times A_1$ (i.e. all zeros outside the
diagonal).
\end{example}

The polynomial Cartan matrices should be regarded as a refinement of the notion
of generalized Cartan matrices in the following sense.
\begin{example} Let $C$ be any generalized Cartan matrix. 
To $C$ we can associate the matrix $P=(p_{ij})$ given
by $p_{ij}(x)=(x-1)^{1-a_{ij}}$ for $i\neq j$ and anything on the diagonal.
Then $P$ is a polynomial Cartan matrix and $C(P)=C$.
More generally, we could take $p_{ij}(x)=(x-q^{a_{ij}})(x-q^{a_{ij}+2})\cdots (x-q^{-a_{ij}})$,
where $q\in\K\backslash\{0\}$.
Assuming that we had a locally finite TGWA, with $R$ a domain, whose polynomial
Cartan matrix was equal to such a $P$, and if all $\mu_{ij}=1$, relations
\eqref{eq:genserre3},\eqref{eq:genserre3Y} would
become the usual (quantum) Serre relations occuring in the (quantum) enveloping
algebra of the Kac-Moody algebra associated to $C$.
In the next section we construct such an algebra
in the case when $C$ is symmetric.
\end{example}

It would be interesting to find conditions under which a locally finite
TGWA $A$ is determined up to isomorphism by its polynomial Cartan matrix $P_A$.
Part of this question would be to determine under what conditions
the generalized Serre relations \eqref{eq:genserre3},\eqref{eq:genserre3Y}
generate the ideal $I$.
In Section \ref{sec:caseA2} we give some sufficient conditions for this
when $n=2$.

\section{A construction of locally finite TGWAs with specified polynomial Cartan matrix}
\label{sec:cartan}

Let $C=(a_{ij})_{i,j=1}^n$ be a symmetric generalized Cartan matrix,
let $q\in\K\backslash\{0\}$ and let $\mu=(\mu_{ij})_{i,j=1}^n$
be a symmetric $n\times n$ matrix without diagonal with $\mu_{ij}\in\K\backslash\{0\}$.

Take $R$ to be the following polynomial algebra over $\K$:
\[R=\K[H_{ij}^{(k)}\mid 1\le i<j\le n, \text{ and } k=a_{ij}, a_{ij}+2, \ldots, -a_{ij}].\]
Define $\si_1,\ldots,\si_n\in\Aut_\K(R)$ by setting, for all $i<j$
and $k=a_{ij}, a_{ij}+2, \ldots, -a_{ij}$:
\begin{align}
\label{eq:cartantgwa_sigmadef1}
\si_j(H_{ij}^{(k)})&=\mu_{ij} q^k H_{ij}^{(k)}+H_{ij}^{(k-2)}, \qquad
 \text{where $H_{ij}^{(a_{ij}-2)}:=0$},\\
\label{eq:cartantgwa_sigmadef2}
\si_i(H_{ij}^{(k)})&=\mu_{ij}^2\si_j^{-1}(H_{ij}^{(k)}),\\
\si_r(H_{ij}^{(k)})&=H_{ij}^{(k)},\qquad r\neq i,j.
\end{align}
For notational purposes, put
\[H_{ij}=H_{ij}^{(-a_{ij})},\quad
 H_{ji}=\si_j^{-1}(H_{ij}) \quad\text{for all $i<j$,
 \quad and $H_{ii}=1$ for all $i$}\]
and define
\[t_i=H_{i1}H_{i2}\cdots H_{in}, \quad \text{for $i=1,\ldots,n$.}\]

\begin{prpdef}
The data $(R,\si,t,\mu)$ satisfies the conditions required in the
definition of a TGWA, namely, the $\si_i$ commute with eachother, and
the consistency relation \eqref{eq:consistency} hold.
We let $\mathcal{T}_{q,\mu}'(C)=A'(R,\si,t,\mu)$
denote the associated twisted generalized Weyl
construction and $\mathcal{T}_{q,\mu}(C)=A(R,\si,t,\mu)$
the corresponding twisted generalized Weyl algebra.
\end{prpdef}
\begin{proof}
It is easy to check that the automorphisms $\si_i$ commute, since on
each $H_{ij}^{(k)}$ either one of the automorphisms is the identity,
or their composition is a multiple of the identity.
To prove the consistency relations
\[t_it_j=\mu_{ij}^2\si_j^{-1}(t_i)\si_i^{-1}(t_j) \quad\forall i\neq j,\]
we can assume by symmetry that $i<j$. Then the right hand side equals
\[\mu_{ij}^2H_{i1}\cdots \si_j^{-1}(H_{ij})\cdots H_{in}
\cdot H_{j1}\cdots \si_i^{-1}(H_{ji})\cdots H_{jn}\]
which equals $t_it_j$ since $\si_j^{-1}(H_{ij})=H_{ji}$ by
choice of notation and
$\mu_{ij}^2\si_i^{-1}(H_{ji})=\mu_{ij}^2\si_i^{-1}\si_j^{-1}(H_{ij})= H_{ij}$
by \eqref{eq:cartantgwa_sigmadef2}.
\end{proof}

To formulate the next theorem, recall that for any $q\in\K\backslash\{0\}$, the
$q$-binomial coefficients
$\begin{bmatrix} m \\ k\end{bmatrix}_q$
may be defined as the elements in the subring $\Z[q,q^{-1}]$ of
$\K$ given by requiring the identity
\begin{equation}\label{eq:qbinomialsdef}
(x+q^{-m+1})(x+q^{-m+3})\cdots (x+q^{m-1})
 =\sum_{k=0}^m \begin{bmatrix} m \\ k\end{bmatrix}_q x^k 
\end{equation}
in $\K[x]$ for any $m\in\Z_{>0}$,
 and $\begin{bmatrix} 0 \\ 0\end{bmatrix}_q =1$ by convention.
If $q$ is not a root of unity and $\K$ has characteristic zero, we have
\begin{equation}\label{eq:qbinomials}
\begin{bmatrix}n\\k\end{bmatrix}_q = \frac{[n]_q!}{[n-k]_q![k]_q!},
\quad [n]_q!=[n]_q[n-1]_q\cdots [1]_q,
\quad [n]_q=q^{-n+1}+q^{-n+2}+\cdots+q^{n-1},
\end{equation}
which can be proved by verifying that both definitions
solve the same recursion relation
\[\begin{bmatrix}m\\k\end{bmatrix}_q =
q^{-k}
\begin{bmatrix}m-1\\k\end{bmatrix}_q +
q^{m-k}
\begin{bmatrix}m-1\\k-1\end{bmatrix}_q.\]

We now come to the main theorem in this section,
which relates the algebra $\mathcal{T}_{q,\mu}(C)$ to the quantum
group $U_q(\mathfrak{g}(C))$ associated to the Kac-Moody algebra $\mathfrak{g}(C)$
(see for example \cite{KS} for an introduction to Drinfel'd-Jimbo quantum groups).

\begin{thm}\label{thm:TqmuC}
\begin{enumerate}
\item[a)] $\mathcal{T}_{q,\mu}(C)$ is locally finite.
\item[b)] The off-diagonal entires of the polynomial Cartan matrix $P=(p_{ij})$
 of $\mathcal{T}_{q,\mu}(C)$ are given by
\begin{align}\label{eq:Tminpol}
p_{ij}(x) &=(x-\mu_{ij}q^{a_{ij}})(x-\mu_{ij}q^{a_{ij}+2})\cdots
(x-\mu_{ij}q^{-a_{ij}})=\\
\label{eq:Tthmid}
&=\sum_{k=0}^{1-a_{ij}}(-1)^{k}\mu_{ij}^k
 \begin{bmatrix}1-a_{ij}\\ k\end{bmatrix}_q x^{1-a_{ij}-k}
\end{align}
for any $i\neq j$.
\item[c)] In $\mathcal{T}_{q,\mu}(C)$, the
quantum Serre relations of the
quantum group $U_q\big(\mathfrak{g}(C)\big)$ hold. That is,
for any $i\neq j$,
\begin{align}
\label{eq:Serreelt1}
&\sum_{k=0}^{1-a_{ij}} (-1)^k 
 \begin{bmatrix}1-a_{ij}\\k\end{bmatrix}_q X_i^{1-a_{ij}-k} X_j X_i^k = 0,\\
\label{eq:Serreelt2}
&\sum_{k=0}^{1-a_{ij}} (-1)^k    
 \begin{bmatrix}1-a_{ij}\\k\end{bmatrix}_q Y_i^k Y_j Y_i^{1-a_{ij}-k} = 0.
\end{align}
Moreover, for any $i\neq j$ and any $m<1-a_{ij}$
the elements $X_i^{m-k}X_jX_i^k$ and $Y_i^{m-k}Y_jY_i^k$,
 $k=0,1,\ldots,m$, are linearly independent over $\K$.
\end{enumerate}
\end{thm}
\begin{proof}
a) Giving $R$ the natural $\Z_{\ge 0}$-gradation by letting all $H_{ij}^{(k)}$
be of degree one, we observe that all automorphisms $\si_i$ preserve the
degree one subspace of $R$. Since the $t_i$ have degree $n$, each $V_{ij}$
is therefore contained in the finite-dimensional subspace of $R$ of elements
of degree $n$, so $\mathcal{T}_{q,\mu}(C)$ is locally finite.

b) Let $1\le i<j\le n$, and consider the following linear subspace of $R$:
\[W_{ij}=\K H_{ij}^{(a_{ij})}\oplus \K H_{ij}^{(a_{ij}+2)}
 \oplus \cdots \oplus \K H_{ij}^{(-a_{ij})}.\]
It has dimension $1-a_{ij}$.
By \eqref{eq:cartantgwa_sigmadef1} the $\K$-algebra automorphism $\si_j$
preserves the subspace $W_{ij}$ and the matrix of $\si_j|_{W_{ij}}$ in
the ordered basis $(H_{ij}^{(a_{ij})}, H_{ij}^{(a_{ij}+2)},\cdots,
  H_{ij}^{(-a_{ij})})$ is given by
\[
\begin{bmatrix}
\mu_{ij}q^{a_{ij}} &            1                 &              &  &  \\
                        & \mu_{ij}q^{a_{ij}+2}    &    1         &  &  \\
         & & \ddots &  \\
         & & & \mu_{ij}q^{-a_{ij}-2} & 1 \\
         & & & & \mu_{ij}q^{-a_{ij}}
\end{bmatrix}
\]
(zeros omitted).
Thus the minimal polynomial
$p_{ij}(x)$ of $\si_j|_{W_{ij}}$ equals
\[
p_{ij}(x)=(x-\mu_{ij}q^{a_{ij}})(x-\mu_{ij}q^{a_{ij}+2})\cdots
(x-\mu_{ij}q^{-a_{ij}}).\]
Since $\si_j(H_{rs}^{(k)})=H_{rs}^{(k)}$ if $r,s\neq j$ it follows that
$p_{ij}$ is also the minimum polynomial for $\si_j$ acting on
the space $V_{ji}$ spanned by all $\si_j^k(t_i)$ ($k\in\Z$).
On $W_{ij}$ we also have, using \eqref{eq:cartantgwa_sigmadef2},
\begin{align*}
&p_{ij}(\si_i)|_{W_{ij}} =
(\mu_{ij}^2\si_j^{-1}-\mu_{ij}q^{a_{ij}})
(\mu_{ij}^2\si_j^{-1}-\mu_{ij}q^{a_{ij}+2})
\cdots
(\mu_{ij}^2\si_j^{-1}-\mu_{ij}q^{-a_{ij}}) =\\
&=\mu_{ij}^{1-a_{ij}}\si_j^{a_{ij}-1}
(\mu_{ij}q^{-a_{ij}}-\si_j)
(\mu_{ij}q^{-a_{ij}-2}-\si_j)
\cdots
(\mu_{ij}q^{a_{ij}}-\si_j) =\\
&=(-\mu_{ij})^{1-a_{ij}}\si_j^{a_{ij}-1} \big(p_{ij}(\si_j)\big)=0.
\end{align*}
On the other hand, $\si_i|_{W_{ij}}$ cannot satisfy any polynomial
of lower degree since it is a multiple of the inverse of $\si_j|_{W_{ij}}$
This shows that $p_{ij}(x)$ is also the minimal polynomial of
$\si_i|_{W_{ij}}$ and thus the minimal polynomial for $\si_i|_{V_{ij}}$.
This proves that the polynomial Cartan matrix $P$ of $\mathcal{T}_{q,\mu}(C)$
is given by \eqref{eq:Tminpol}.
The equality \eqref{eq:Tthmid} follows from the definition \eqref{eq:qbinomialsdef}
of the $q$-binomial coefficients.

c) This follows from part b) and Theorem \ref{thm:lfthm}.
\end{proof}

We believe that the quantum Serre relations
\eqref{eq:Serreelt1},\eqref{eq:Serreelt2}
in fact generate the ideal $I$ of the TGWC $\mathcal{T}_{q,\mu}'(C)$.
This would yield a complete presentation of $\mathcal{T}_{q,\mu}(C)$
by generators and relations.
The result in the next section implies that this is true in the case
when $C$ is the Cartan matrix of type $A_2$.

\section{Sufficient conditions for Serre-type relations to generate the ideal $I$ of a TGWC}
\label{sec:caseA2}

\begin{thm}\label{thm:suffcond}
Assume that $A'=A'(R,\si,t,\mu)$ is a
twisted generalized Weyl construction of rank $n=2$
where $R$ is a domain, with the following properties:
\begin{enumerate}
\item[(P1)] For any positive integer $k$ and any $r\in R$ the following
implication hold:
\[r
\cdot  (\si_1\si_2)^k(t_1)\cdot (\si_1\si_2)^{k-1}(t_1)\cdots 
 (\si_1\si_2) (t_1)
\in R\cdot \si_1(t_1) \Longrightarrow r\in R\cdot \si_1(t_1),\]
\item[(P2)] the ideal $I$ of $A'$ contains elements of the form
\begin{align*}
s_1 &:= X_2X_1^2 - \xi_1  X_1X_2X_1 - \xi_2 X_1^2X_2,\\
s_2 &:= X_2^2X_1 - \eta_1 X_2X_1X_2 - \eta_2 X_1X_2^2,
\end{align*}
where $\xi_i,\eta_i$ are nonzero elements of $\K$.
\end{enumerate}
Then $I$ is generated by $s_1, s_2, s_1^\ast, s_2^\ast$.
\end{thm}

\begin{rem}
Condition (P1) is satisfied if the following two properties hold:
\begin{enumerate}
\item[(P1a)] $\si_1\si_2(t_1)=\la t_1$ for some $\la\in\K\backslash\{0\}$,
and 
\item[(P1b)] $R\cdot t_1$ is a prime ideal in $R$ and the ideals
$R\cdot t_1$ and $R\cdot \si_1(t_1)$ are coprime, i.e. their
sum is all of $R$.
\end{enumerate}
\end{rem}

\begin{proof}
We begin with a simple observation.
Let $g=(g_1,g_2)\in \Z^2$ and suppose $a\in I\cap A_g$, $a\neq 0$.
If $g_1\ge 0$ and $g_2\le 0$ then $a=r X_1^{g_1} Y_2^{g_2}$ for some
$r\in R$. Since $I$ is completely gr-prime and contains no monic
monomials (Proposition \ref{prp:completely_gr-prime} and
Corollary \ref{cor:nomonicmonomials})
we get $r\in I$. But $I\cap R=0$ so $r=0$, contradicting $a\neq 0$.
Similarly, we cannot
have $g_1\le 0$ and $g_2\ge 0$. Thus, if $I\cap A_g\neq 0$, then either
$g_1>0$ and $g_2>0$, or $g_1<0$ and $g_2<0$.

Now let $I_S$ denote the ideal in $A'$ generated by $s_1,s_2,s_1^\ast,s_2^\ast$.
By Property (P2) we know $I_S\subseteq I$.
Assume that $I_S\subsetneq I$. Then, among all homogenous
elements in $I\backslash I_S$,  choose one, $a$, with minimal
total degree $|g_1+g_2|$, where $(g_1,g_2)=\deg a \in \Z^2$.
After applying the involution $\ast$ if necessary, we can,
by the previous paragraph, assume that $g_1,g_2>0$.

For a monic monomial
$w=Z_1\cdots Z_k\in A'$ where $Z_i\in\{X_1, X_2\}$ $\forall i$, define the \emph{length}
$\ell(w)$ as the number of pairs $(i,j)$, $1\le i<j\le k$ such that
$Z_i=X_2$ and $Z_j=X_1$. The element $a$ can be written as a sum of monomials
in the noncommuting elements $X_1, X_2$ with coefficients from
$R$ written on the left. Some of these monomials can be reduced
mod $I_S$ with the reductions
\begin{align}
X_2X_1^2 &\rightarrow \xi_1  X_1X_2X_1+ \xi_2 X_1^2X_2,\\
X_2^2X_1 &\rightarrow \eta_1 X_2X_1X_2+ \eta_2 X_1X_2^2.
\end{align}
At each reduction a monomial $w$ is replaced by a sum of two monomials
each having strictly lower length than $w$.
Thus after finite number of steps we obtain an element of the
form
\begin{equation}\label{eq:aelement}
a'=\sum_{i=0}^{\min(g_1,g_2)} \be_i X_1^{g_1-i}(X_2X_1)^i X_2^{g_2-i}
\end{equation}
where $\be_i\in R$. We can without loss of generality assume
$a=a'$ since we only added elements from $I_S$ of degree $(g_1,g_2)$.
Then we have
\[
aY_1^{g_1} = \sum_{i=0}^{\min(g_1,g_2)}
 \al_i X_1^{g_1-i} \cdot (\si_1\si_2)^i(t_1)\cdot (\si_1\si_2)^{i-1}\cdots
  \si_1\si_2(t_1)\cdot Y_1^{g_1-i}X_2^{g_2},
\]
where $\al_i$ equals $\be_i$ times some power of $\mu_{21}$.
Multiplying further by $Y_2^{g_2}$ from the right we get zero
since $a\in I$. After that we can cancel $X_2^{g_2}Y_2^{g_2}$
since it is nonzero and $R$ has no zero-divisors. Thus
\begin{equation}\label{eq:sum}
\sum_{i=0}^{\min(g_1,g_2)}
 \al_i \si_1^{g_1}\si_2^i(t_1)\cdot \si_1^{g_1-1}\si_2^{i-1}(t_1)\cdots 
 \si_1^{g_1-i+1}\si_2 (t_1)\cdot 
 X_1^{g_1-i}Y_1^{g_1-i}=0.
\end{equation}
For simplicity we now assume $g_1\le g_2$. The case $g_1>g_2$
can be treated similarly.
Then observe that, for $i<g_1$, the $i$-term in the sum
\eqref{eq:sum} contains a factor $\si_1(t_1)$, namely by
simplifying $X_1^{g_1-i}Y_1^{g_1-i}$.
So
\[\al_{g_1}\cdot  (\si_1\si_2)^{g_1}(t_1)\cdot (\si_1\si_2)^{g_1-1}(t_1)\cdots 
 (\si_1\si_2) (t_1)
  + \ga \cdot \si_1(t_1) =0\]
for some $\ga\in R$. By Property (P1) this implies that
 $\al_{g_1}=\tilde \al_{g_1}\cdot \si_1(t_1)$
for some $\tilde \al_{g_1}\in R$.
Substituting this in the last term of \eqref{eq:aelement} 
and using that $\si_1(t_1)X_2X_1 =
\mu_{21}^{-1}\si_1^{-1}\si_2^{-1}(t_1)X_1X_2$
we get
\[a=\sum_{i=0}^{g_1-1} \al_iX_1^{g_1-i}(X_2X_1)^iX_2^{g_2-i}
+\tilde\al_{g_1}
\mu_{21}^{-1}\si_1^{-1}\si_2^{-1}(t_1)
 X_1X_2(X_2X_1)^{g_1-1}X_2^{g_2-g_1} = X_1 b\]
for a certain $b\in A'$. 
Since $X_1\notin I$, we conclude by Proposition \ref{prp:completely_gr-prime}
that $b\in I$. But then $b$ has total degree $g_1+g_2-1$.
Thus $b\in I_S$, otherwise 
it contradicts the choice of $a$ as an element
of $I\backslash I_S$ with minimal total degree.
But then $a=X_1b\in I_S$ also which is a contradiction.
This shows that $I_S\subsetneq I$ is false and thus $I_S=I$.
\end{proof}

\begin{example}[The TGWA of ``type $A_2$'', contd.]
\label{ex:typeA2generators}
Recall from Example \ref{ex:typeA2} that
 $n=2$, $R=\K[H]$, $\si_1(H)=H+1$, $\si_2(H)=H-1$,
 $t_1=H$, $t_2=H+1$, $\mu_{12}=\mu_{21}=1$ and
let $A=A(R,\si,t,\mu)$ be the TGWA. Conditions (P1a), (P1b) and (P2)
above are satisfied and thus $A$ is isomorphic to the $\K$-algebra
with generators $X_1,X_2,Y_1,Y_2,H$ and defining relations
\begin{gather*}
\begin{alignedat}{2}
X_1H&=(H+1)X_1, &\qquad  X_2H&=(H-1)X_2, \\
Y_1H&=(H-1)Y_1, &\qquad Y_2H&=(H+1)Y_2, \\
Y_1X_1&=X_2Y_2=H, &\qquad   Y_2X_2&=X_1Y_1=H+1,
\end{alignedat}\qquad
\begin{aligned}
X_1^2X_2-2X_1X_2X_1+X_2X_1^2&=0,\\
Y_1^2Y_2-2Y_1Y_2Y_1+Y_2Y_1^2&=0.
\end{aligned}
\end{gather*}
It is interesting to notice that $A$ contains
two copies of the first Weyl algebra (with generators $X_1,Y_1$ and $X_2,Y_2$
respectively), and at the same time $X_1,X_2$ (and $Y_1,Y_2$) generate
a subalgebra isomorphic to the enveloping algebra of the positive (or negative)
nilpotent part of $\mathfrak{sl}_3(\K)$.
\end{example}

\begin{cor}\label{cor:A2}
Let $C$ be the Cartan matrix of type $A_2$.
Then the maximal graded ideal with trivial zero component, $I$,
of $\mathcal{T}_{q,\mu}'(C)$ is generated by
the elements
\begin{align*}
X_1^2X_2-(q+q^{-1})X_1X_2X_1+X_2X_1^2,\\
Y_1^2Y_2-(q+q^{-1})Y_1Y_2Y_1+Y_2Y_1^2.
\end{align*}
\end{cor}
\begin{proof}
Properties (P1a) and (P1b) are easily verified for the TGWC
$\mathcal{T}_{q,\mu}'(C)$, and Property (P2) holds by
Theorem \ref{thm:TqmuC}c), so Theorem \ref{thm:suffcond} can
be applied and yields the result.
\end{proof}
It is now straightforward to write down a set of generators and relations
for $\mathcal{T}_{q,\mu}'(C)$ when $C$ is of type $A_2$.

\noindent\textsc{
University of Amsterdam,
Korteweg-de Vries Institute for Mathematics,
P.O. Box 94248,
1090 GE Amsterdam,
The Netherlands}

\noindent\texttt{Email: jonas.hartwig@gmail.com}

\end{document}